\documentclass[a4paper,10pt]{article}
\def\uu{\bigsqcup}
\def\eps{\varepsilon}
\def\Z{\bf Z}
\def\N{\bf N}
\usepackage[T2A]{fontenc}
\usepackage[cp1251]{inputenc}
\usepackage[english,russian]{babel}
\usepackage[tbtags]{amsmath}
\usepackage{amsfonts,amssymb}

\usepackage{mathrsfs}

\usepackage{graphicx}

\date{}
\begin{document}
\Large

\title{ Спектральные свойства \\ и приближение присоединений    \\ 
бесконечных действий  ранга один
}
\author{ Рыжиков В.В.}

\maketitle

\begin{abstract}
    В заметке доказано, что эргодическое самоприсоединение бесконечного преобразования ранга 1
является частью  слабого предела сдвигов диагональной меры. Предложен континуальный класс неизоморфных 
преобразований с полиномиальным замыканием. Эти преобразования обладают 
минимальными самоприсоединениями и некоторыми необычными спектральными свойствами. Так, например,
тензорные произведения степеней преобразования имеют как сингулярный так и 
лебеговский спектр в зависимости от выбора степеней.

\end{abstract}

\section{ Аппроксимации самоприсоединений преобразования ранга 1}
Действия ранга 1 в эргодической теории представляют интерес по разным причинам.
Они позволяют строить  динамические системы с разнообразными алгебраическими,
асимптотическими и спектральными свойствами (см., например, \cite{Ru},\cite{Ry}).
 Cтруктуры самоприсоединений (self-joinings), факторов и централизатора действия ранга 1  связаны с  
 его слабым замыканием.   Например, в
\cite{K} доказано, что  для пространства с конечной мерой централизатор преобразования ранга один  
лежит в слабом замыкании его степеней.    
Аналог этого утверждения  в случае бесконечного пространства был установлен только 
для перемешивающих действий   \cite{RT} (см. также \cite{KR}).  

В предлагаемой заметке   будет получена частичная аппроксимируемость  
эргодических    самоприсоединений  бесконечного преобразования ранга 1 сдвигами диагональной меры. 
Этот факт  обобщает результат работы \cite{RT} и  
применяется  к классу неперемешивающих конструкций с полиномиальным слабым замыканием: 
для  них получено свойство 
минимальных самоприсоединений. Предъявляется континуум  неизоморфных преобразований  
с  полиномиальным слабым замыканием.
Они обладают  сингулярным спектром, но
 попарные свертки их спектральных мер являются лебеговскими. Гауссовские и пуассоновские надстройки над ними 
также обладают   сингулярным спектром, а  тензорные произведения  надстроек имеют 
лебеговскую компоненту в спектре.  Модификация
конструкций позволяет получать разнообразие спектральных свойств тензорных произведений 
степеней  преобразования.

Напомним необходимые определения.  Сохраняющее меру обратимое преобразование
$ T: X \to X $ пространства Лебега $ (X, \mu) $ имеет {\it ранг} 1, 
если существует последовательность 
$$ \xi_j= \{E_j, \ TE_j, \ T^2 E_j, \ \dots, T^{ h_j} E_j, \tilde {E}_j \} $$ 
измеримых разбиений пространства 
$ X $ таких, что 
$ \xi_j $ стремится к разбиению на точки.  Это означает, что любое множество 
конечной меры сколь угодно точно
для всех достаточно больших $j$  приближается $ \xi_j$-измеримыми  множествами.

  {\it Cамоприсоединением} преобразования $T$,  называется мера $\nu$ на 
$X\times X$, для которой выполнены следующие свойства:
\begin{enumerate}
  \item $ (T\times T)\nu = \nu$,
  \item $\nu(X\times A) = \nu(A\times X) = \mu(A)$ для любого множества $A$ конечной меры.
\end{enumerate}
Пусть обратимое сохраняющее меру преобразование $S$ коммутирует с $T$. Рассмотрим меру 
$\Delta_S$ на $(X\times X)$, определенную
равенством $\Delta_S=(Id\times S)\Delta$, где $\Delta$ -- диагональная мера: 
$\Delta(A\times B)=\mu(A\cap B)$. Так как $S$ коммутирует с $T$, то мера 
$\nu=\Delta_S$ удовлетворяет условию инвариантности $(T\times T)\nu=\nu$. 

Отображение $x \rightarrow (x,Sx), x \in X$, задает изоморфизм между 
системами $(T,X,\mu)$ и $(T\times T, X\times X, \Delta_S)$. 
Следовательно, если мера $\mu$ эргодическая, то мера $\Delta_S$ также будет эргодической. 
В частности, эргодическими самоприсоединениями эргодического преобразования $T$  будут 
все меры вида $$\Delta^n=\Delta_{T^n}.$$

Если же у преобразования нет  эргодических  самоприсоединений, отличных от $\Delta_{T^n}$,  говорят, 
что оно обладает  {\it минимальными самоприсоединениями (minimal self-joinings)}. 
Действие с минимальными самоприсоединениями  обладает тривиальным централизатором, то есть коммутирует
только со своими степенями.
\vspace{3mm}

\bf Теорема 1. \it
Если $\nu$ -- эргодическое самоприсоединение бесконечного преобразования $T$  ранга 1, 
то для некоторой 
последовательности $k(j)$ и некоторой меры $\nu'$  выполнено
$$\Delta^{k(j)} \to \frac 1 2 \nu +\nu'.$$ \rm

\vspace{3mm}

Оператор $P$  в $L_2(X,\mu)$ называем марковским, если $P$ положителен, т.е. сопоставляет 
неотрицательным функциям  неотрицательные,  и операторы  $P, P^\ast$ сохраняют интеграл. 
Самоприсоединениям 
преобразования $T$  отвечают марковские операторы, 
коммутирующие с $T$ (в заметке $T$ обозначает и преобразование и соответствующий ему 
оператор:  $Tf(x)=f(Tx)$). Эта связь задается формулой 
$$  (Pf,g)=\int_{X\times X}  f\otimes g d\nu. $$ 

На операторном языке утверждение теоремы 1 имеет вид 
$$T^{k(j)} \to \frac 1 2 P +P',$$
где марковский  оператор $P$ отвечает  самоприсоединению $\nu$,  а $P'$ --  некоторый положительный оператор.

\bf Следствие. \it Если  слабое замыканине действия   преобразования $T$ ранга 1 полиномиально, т.е. 
состоит из линейных комбинаций степеней этого преобразования, то $T$ обладает минимальными 
самоприсоединениями. \rm
\vspace{3mm}

Действительно, эргодическое самоприсоединение $\nu$  есть 
некоторая линейная комбинация 
мер  $\Delta^k$,$k\in\Z$, но в силу его эргодичности оно может состоять только из одной 
эргодической компоненты.
Следовательно,  $\nu=\Delta^n$.

\section {Доказательство теоремы 1} Так как $T$ -- преобразование ранга один,
 найдется последовательность разбиений $ \{E_j, \ TE_j, \ T^2 E_j, \ \dots, T^{ h_j} E_j, \tilde {E}_j \}$, 
стремящаяся к разбиению на точки.  
Определим множества 
$$ C^k_j=\uu_{i=0}^{h_j-k} T^{i+k} E_j\times T^{i} E_j, \ \ 
k\in[0,h_j],$$
$$ C^k_j=\uu_{i=0}^{h_j+k}T^{i} E_j\times T^{i-k} E_j, \ \ k\in[-1,-h_j].$$
Обозначим  $\Delta^k=\Delta_{T^k}$, 
определим меру $\Delta^k_j$  равенством
$$\Delta^k_j (A\times B)= \Delta^k ((A\times B)\cap C^k_j),\ \ k\in[-h_j,h_j],
$$
где $A,B$ -- множества конечной меры.
В работе \cite{RT} показано, что
   $$
\sum_{|k|\leq h_j} a_j^k\Delta^k_j(A\times B) \to \nu(A\times B), \ \ j\to\infty, $$
причем будут  выполнены условия
$$
a_j^k \geq 0, \ \ \sum_{|k|\leq h_j} a_j^k\to a,   \  1\leq a \leq 2.$$
Положив   $c_j^k=a_j^k/a$,
перепишем приведенное  утверждение в  виде
   $$
\sum_{|k|\leq h_j} c_j^k\Delta^k_j(A\times B) \to c\nu(A\times B), \ \ j\to\infty, $$
где
$$\ \ c_j^k\geq 0, \ \ \sum_{|k|\leq h_j} c_j^k=1.   $$

Пусть далее   последовательность $c_j^k$, удовлетворяющая указанным условиям,  выбирана так, чтобы 
 число $c$  было максимальным ($\frac 1 2 \leq c\leq 1$). Положим
$$ D(\eps,j):=\{k:\Delta^k_j(A\times B) - 
c\nu(A\times B)>\eps\}.$$
Покажем, что для  множеств $A,B$ конечной меры, $\nu(A\times B)>0$,  выполняется
 $$\sum_{k\in  D(\eps,j) } c_j^k\to 0. $$
Если это не так, то для некоторой подпоследовательности  $j_n$,  которую снова обозначим через   $j$,
имеем $$\sum_{k\in  D(\eps,j) } c_j^k\to c'>0.$$  Из этого, как мы увидим, следует, что 
выбранное число $c$ не является
 максимальным.  Действительно, положив  $\tilde c_j^k=\frac { c_j^k}{c'}$ для $k\in  D(\eps,j)$ 
и $\tilde c_j^k=0$ для остальных $k$,   получим $$\sum_{|k|\leq h_j} \tilde c_j^k=1, \ 
\sum_{k\in  D(\eps,j) } \tilde c_j^k\Delta^k_j\to \tilde c\nu.\ $$
Сходимость  к $c\nu$ имеет место по причине того, что  предельная мера инвариантна относительно 
$T\times T$ и абсолютно непрерывна
относительно эргодической меры $\nu$. При $\nu(A\times B)>0$ имеем 
$$\sum_{k\in  D(\eps,j) } \tilde c_j^k\Delta^k_j(A\times B) \geq  (c+\eps)\nu(A\times B),\ $$
следовательно,  $\tilde c > c$. Это противоречит максимальности числа $c$.

Обозначая $$
\bar D(\eps,j):=\{k:\Delta^k_j(A\times B) -  c\nu(A\times B)<-\eps\},$$
получим 
$$\sum_{k\in  \bar D(\eps,j) } c_j^k\to 0. 
$$
Иначе некоторая подпоследовательность этих сумм болше, чем $c''>0$.  Тогда,   
переходя к рассмотрению комбинаций мер  $\Delta^k_j$ для номеров $k$ вне множества 
$\bar D(\eps,j)$, получим, что в среднем для них выполняется
$$
\Delta^k_j(A\times B) -  c\nu(A\times B)\geq c''\eps .$$
Это  противоречит  максимальности числа $c$.

Таким образом, для весового большинства номеров $k$  (когда сумма $c^k_j$ по таким $k$
близка к  1) значения $\Delta^k_j(A\times B)$ близки к   $c\nu(A\times B)$.
Диагональным методом для любого счетного набора
пар  $A,B$, множеств конечной меры, удовлетворяющих условию $\nu(A\times B)>0$,
 можно выбрать последовательность $k(j)$ такую, что  
   $$
\Delta^{k(j)}_j(A\times B) \to c\nu(A\times B), \ c\geq \frac 1 2 .
$$
Пусть такой набор плотен в семействе всех пар  множеств конечной меры.
В силу проекционных свойств рассматриваемых мер  сходимость имеет место  для 
любых  пар множеств  $A,B$ конечной меры.   Отсюда сразу следует, что 
$$\Delta^{k(j)} \to \frac 1 2 \nu +\nu'.$$ 
Теорема доказана.
   
\section{Действия   с полиномиальным  слабым замыканием}
 Бесконечные преобразования $T$, обладающие  свойством
$T^k\to_w 0$ при $k\to\infty$, называются перемешивающими. Для них
полугруппа  слабых пределов степеней $Lim(T^n)$  есть $\{0\}.$ 
 Мы рассмотрим  преобразования 
с нетривиальной  полиномиальной полугруппой.  Пример такого преобразования $T$ 
 имеется  в  заметке \cite{KuR}, где доказано, что   $$Lim(T^n)=
\{2^{-m}T^s:  m\geq 1, s\in \Z  \}\cup \{0\}.$$
В этом параграфе будет приведено  континуальное семейство попарно неизоморфных 
преобразований с  полиномиальным слабым замыканием.  
\vspace{3mm}

\bf Конструкция преобразования ранга 1. \rm 
 Фиксируем натуральное число $h_1$ (высота башни на этапе 1), последовательность $r_j\to\infty$ 
($r_j\geq 2$ -- число колонн, на которое разрезается  башня на этапе $j$)
и последовательность  целочисленных наборов 
$$ \bar s_j=(s_j(1), s_j(2),\dots, s_j(r_j-1),s_j(r_j)), \ s_j(i)\geq 0.$$ 
Параметры $r_j$, $\bar s_j$ полностью  определяют конструкцию преобразования ранга 1, 
напомним ее описание.

Пусть на  шаге $j\geq 1$ определена 
система  непересекающихся полуинтервалов, одинаковой длины (башня высоты $h_j$)
$$E_j, TE_j, T^2E_j,\dots, T^{ h_j-1}E_j,$$
причем на на полуинтервалах  $E_j, TE_j, \dots, T^{ h_j-2}E_j$
пребразование $T$ действует как   перенос полуинтервалов.  

 Представим 
 $E_j$ как дизъюнктное объединение  полуинтервалов $ E_j^i$, $1\leq i \leq r_j$, одинаковой длины. 
Набор  
$E_j^i, TE_j^i ,T^2 E_j^i,\dots, T^{h_j-1}E_j^i$ называется 
$i$-ой колонной  этапа $j$.
К этому набору добавим $s_j(i)$ полуинтервалов меры $\mu(E_j^i)$, 
 получая набор непересекающихся полуинтервалов
$$E_j^i, TE_j^i, T^2 E_j^i,\dots, T^{h_j-1}E_j^i, T^{h_j}E_j^i, T^{h_j+1}E_j^i, \dots, T^{h_j+s_j(i)-1}E_j^i.$$

  Терерь соберем   надстроенные колонны в одну башню, для этого положим 
$$T^{h_j+s_j(i)}E_j^i = E_j^{i+1}$$ при $i<r_j.$
Обозначая  $E_{j+1}= E^1_j$, получили  башню этапа $j+1$:
$$E_{j+1}, TE_{j+1}, T^2 E_{j+1},\dots, T^{h_{j+1}-1}E_{j+1},$$
где 
$$ h_{j+1} =h_jr_j +\sum_{i=1}^{r_j}s_j(i).$$

Продолжая построение до бесконечности,  определяем  преобразование $T$  
 на объединении  $X$ всех рассматриваемых полуинтервалов. 
Преобразование $T$ обратимо и  сохраняет меру Лебега. 
Мера пространства $Х$ бесконечна, если выполнено 
$$\sum_j \sum_{i=1}^{r_j}\frac {s_j(i)}{h_jr_j }=\infty.$$
Преобразование ранга 1 является эргодическим и имеет простой спектр.

{\bf Примеры.} 
Рассмотрим следующий класс конструкций.  Пусть $r_j=r\geq 2$ и
выполнено
$$s_j(1)=0, h_j<< s_j(2)<<\dots<< s_j(r-1)<<s_j(r) \eqno (\ast),$$
$...<<...$ означает,  что отношение левой части к правой стремится к 0 при $j\to\infty$. 
Если последовательность $n(k)$  имеет вид 
$$n(k)=  \pm  h_{j_1(k)} \pm  h_{j_2(k)}\dots \pm  h_{j_p(k)}+s,$$
причем   ${j_1(k)}>  {j_2(k)}>\dots \ >{j_p(k)}\to\infty$, то 
$$T^{n(k)}\to \frac 1 {r^p}T^s.$$
Хотя такие последовательности $n(k)$ при $r>2$  не описывают все последовательности 
с ненулевым пределом для  $T^{n(k)}$,
можно показать, что пределов, отличных от  $0$ и  $\frac 1 {r^p}T^s$, для такой конструкции нет.
Рассмотрим наиболее простой случай. 
\vspace{2mm}

\bf Утверждение 1. \it Пусть конструкция $T$ задана параметрами
$ r_j=2,$   $s_j(1)=0,$  $s_j(2)>>h_j.$  Тогда
 $Lim(T^n)= \{2^{-m}T^s:  m\geq 1, s\in \Z  \}\cup \{0\}.$\rm
\vspace{2mm}

Доказательство.
Фиксируем   $\xi_{j_0}$-измеримые множества $A,B$.
Если  $$\lim_{k\to\infty} \ \mu(T^{n(k)}A\cap B)>0,$$ 
то для достаточно больших $k$ 
для некоторой последовательности $j(k)$ выполнется
 $$h_{j(k)} - 2h_{j(k)-1}< n(k)< h_{j(k)} + 2h_{j(k)-1}.\eqno (3)$$ 
Иначе имеем
$$h_{j(k)} + 2h_{j(k)-1}\leq n(k)\leq h_{j(k)+1 } - 2h_{j(k)}.$$
Но в этом случае получим $\mu(T^{n(k)}A\cap B)=0$, так  как  $A,B$ лежат в  башне этапа $j(k)-1$, 
а эта башня не пересекается с ее образом под действием   $T^{n(k)}$.   

Таким образом,  выполнено (3). 
Обозначим последовательность $j(k)$ через $j_1(k)$.
Заметим, что пересечение
$T^{n(k)}A\cap B$  имеет меру в два раза меньше меры пересечения  $T^{n(k)-h_{j(k)}}A\cap B$. 
Убедимся в этом.  Запишем $A=A_k\uu T^{h_{j_1(k)}}A_k$, где $A_k$ --
пересечение  множества $A$ с первой колонной в башне этапа $j_1(k)$, 
пересечение  множества $A$ со второй колонной есть $T^{h_{j_1(k)}}A_k$.    Аналогично представим 
$B=B_k\uu T^{h_{j_1(k)}}B_k$. Получаем
  $$T^{n(k)}A\cap B=T^{n(k)}A_k\cap T^{h_{j_1(k)}}B_k,$$
$$\mu(T^{n(k)}A\cap B)=\mu (T^{n(k)-h_{j_1(k)}}A_k\cap B_k)= \frac 1 2 \mu(T^{n(k)-h_{j_1(k)}}A\cap B),$$
$$\lim_k\mu(T^{n(k)}A\cap B)=\frac 1 2 \lim_k\mu(T^{n(k)-h_{j_1(k)}}A\cap B).$$
Теперь  в  (3), рассматриваем вместо последовательности  $n(k)$ последовательность 
$n(k)\pm h_{j_1(k)}$.
Если она не ограничена, получим
$$\lim_k\mu(T^{n(k)}A\cap B)=\frac 1 4 \lim_k\mu(T^{n(k) \pm  h_{j_1(k)}\pm  h_{j_2(k)}}A\cap B).$$
Продолжая описанную процедуру, т.е. рассматривая вместо последовательности  $n(k)$ последовательности 
$n(k)\pm h_{j_1(k)}\pm h_{j_2(k)}$ и т.д., для некоторого $m$ получим ограниченную последовательность
 $$n(k) \pm  h_{j_1(k)} \pm  h_{j_2(k)}\dots \pm  h_{j_m(k)}.$$
Иначе для  любых множеств $A,B$ конечной меры 
$$\lim_k\mu(T^{n(k)}A\cap B)\leq  2^{-m(k)}\mu(A) \to 0,$$
что противоречит исходному предположению о том, предел $T^{n(k)}$ ненулевой.

Таким образом, установлено, что  $m(k)$  и  
$s(k):=n(k) \pm  h_{j_1(k)} \pm  h_{j_2(k)}\dots \pm  h_{j_m(k)}$  
являются постоянными для всех больших $k$
и выполняется 
$$T^{n(k)}\to 2^{-m}T^s.$$

\section{Спектральные  свойства конструкций} 
Преобразования, описанные в предыдущем параграфе обладают следующим свойством:
если для некоторых  множеств $A,B$  конечной меры   
$$\lim_k\mu(T^{n(k)}A'\cap B')> 0,\ n(k)\to\infty$$ 
(это верно, когда  $T^{n(k)}$ имеет ненулевой предел),
то для некоторой последовательности ${j(k)}$ имеем $$\frac {n(k)}{h_{j(k)}}\to 1.$$ 
Это вытекает из (3) с учетом того, что $ {h_{j(k)-1}}<<{h_{j(k)}}.$

{ \bf Континуум неизоморфных  преобразований с  полиномиальным слабым замыканием.}
Зафиксируем некоторую конструкцию $T_1$ из утверждения 1.  
Для каждого $a>1$  рассмотрим  конструкцию $T_a$,
у которой соответствующие высоты асимптотически в $a$ раз больше, чем высоты у преобразования $T_1$. 
Рассмотрим $T_a$ и $T_b$  при  $a<b$.
Пусть для некоторого оператора $J$  выполнено
$ T_aJ=JT_b,$  тогда 
$$ T_a^{n}J=JT_b^{n}.$$  Пусть  $h_j$   -- последовательность высот   для $T_b$.  Получим
$$T_a^{h_j}\to 0, \ T_b^{h_j}\to \frac 1 2 I, $$  следовательно, 
 $J=0$.  Тем самым мы установили спектральную дизъюнктность преобразований $T_a$ и $T_b$  при  $a<b$.

Интересен тот факт, что все  произведения $T_a\times T_b$  при  $a<b$  спектрально 
изоморфны между собой.   Напомним, что преобразование  $S$ называется диссипативным,
если все пространство представляется как объединение непересекающихся множеств $S^iD,  i\in \Z$  
для некоторого измеримого  множества $D$. Спектр такого преобразования счетнократный лебеговский.
\vspace{2mm}

\bf Утверждение 2. \it При  $a<b$  произведение $T_a\times T_b$ диссипативно.\rm
\vspace{2mm}

Доказательство.  Для  $\xi_j$-измеримого  множества  $A$,
 $\xi'_j$-измеримого  множества  $A'$ 
(оба конечной меры)   выполнено  
$$(A\times A')\cap  (T_a\times T_b)^k (A\times A')=\phi$$
для   всех   $k>h'_j$ при условии, что $\frac {h_i}{h_i'}$
близко к $\frac a b$  для всех  $i\geq j$. Но это условие выполнено для достаточно больших $j$.
 Таким образом,  $A\times A'$ лежит в диссипативной части.
Увеличивая  $A$ и $A'$,  получим, что диссипативная часть   произведения 
$T_a\times T_b$ совпадает со всем пространством
 $X\times X$.
 
{ \bf Гауссовские и пуассоновские надстройки.}  Если $T$ --   автоморфизм  
стандартного бесконечного пространства 
с мерой,  ему соответствуют  гаусовская $G(T)$ и пуассоновская  $P(T)$  
надстройки (см., например, \cite{N}).   
Автоморфизмы $G(T)$ и  $P(T)$  действуют на вероятностном пространстве, 
как операторы они унитарно эквивалентны  
$$exp(T)=\bigoplus_{n=0}^{\infty} T^{\odot n},$$ где  $T$ --
оператор, отвечающий автоморфизму  $T$, $T^{\odot 0}$ -- одномерный тождественный оператор, 
$T^{\odot n}$ -- симметрическая тензорная степень оператора $T$.

Таким образом, гауссовские и пуассоновские надстройки  обсуждаемых  конструкций $T_a$ 
 имеют чисто сингулярный спектр,
а их попарные тензорные произведения имеют счетнократную лебеговскую компоненту в спектре. 
 Это вытекает из того, что спектр
$T_a^{\odot m}$  сингулярный,  а  спектр   $T_a^{\odot m}\otimes T_b^{\odot n}$ 
 лебеговский  при  $m,n>0$.
\vspace{2mm}

{ \bf  О спектре  произведений  $\bf T\otimes T^n$.}  
Слабая сходимость  $T_a^{h_j}\to \frac 1 2 I $  влечет
за собой сингулярность спектра  (см. \cite{O}, \cite{S}).  Спектр $T_a\times T_a$  также сингулярен,
так как $$(T_a\times T_a)^{h_j}\to \frac I 4 I\otimes I. $$  Но при   $n>1$ произведения 
$T_a\times T_a^n$ 
имеют счетнократный лебеговский 
спектр, так как  они     диссипативны.    Доказательство  диссипативности произведений 
$T_a\times T_a^n$ при   
$n>1$ аналогично доказательству утверждения 2. 

Модификации  конструкций  приводит к  
любопытным  спектральным свойствам, представляющим самостоятельный интерес.  
\vspace{2mm}

\bf  Теорема 2. \it  Для любого $N\geq 2$ найдется  автоморфизма $T$ такой, что  
при  $1\leq m < n\leq N$  спектр произведения $T^m\otimes T^n$     простой и сингулярный, 
тензорные квадраты $T^m\otimes T^m$, $m\neq 0$,  имеют однородный  
сингулярный спектр кратности 2,  а
при $n >N $  спектр   $T\otimes T^n$     счетнократный лебеговский.\rm
\vspace{2mm}

Доказательство.  Подходящая  конструкция ранга 1 задается параметрами:

$r_j=N+1$,

$\bar s_j=(0,0,\dots,0,s_j(N),s_j(N+1))$, 
\\
где  $s_j(N+1)>>h_j$, а последовательнность 
$s_j(N)$ не имеет быстрого роста (скажем, $s_j(N)\leq j$) и каждое положительное значение 
она  принимает бесконечное число раз. Например,  годится последовательность 
$$ 1,2, 1,2,3, 1,2,3,4, 1,2,3,4,5,  \dots.$$ Легко найти континуум незоморфных конструкций, 
удовлетворяющих указанным  требованиям.

Пусть  $1\leq n \leq N$, а $ j'$ пробегает лишь те значения, для которых  $s_{j'}(N)=p$. 
Тогда  выполняется
$$T^{-nh_{j'}}\to  \frac {N-n}{N+1}I+\frac {1}{N+1}T^p. \eqno (4)$$
Вычисление таких пределов стандартно (похожие пределы использовались  в работе \cite{Ry}).
Для отрицательных значений $p$ рассматриваем последовательность $T^{nh_j}$  вместо $T^{-nh_j}$. 
Обозначим через $W_{m,n}$  алгебру фон Неймана, являющуюся 
наименьшей слабо замкнутой алгеброй, порожденной степенями оператора  $T^m\otimes T^n$.
 Из (4)  следует 
$$ ((N-m)I+ T^{p'}) \otimes ((N-n)I+ T^{p'})  \in W_{m,n}\eqno (5)$$
для всех $p'\in \Z$. 
Подставляя в  (4)  $n=N$, получим  пределы $\frac {T^{p'}}{N+1}$, для всех $p'\in \Z$.
Рассматривая замыкание множества  операторов из (5), замечаем, что
$$ \left((N-m)I+ \frac {T^p}{N+1}\right) \otimes \left((N-n)I+ 
\frac {T^p}{N+1}\right)\in W_{m,n}. \eqno (6)$$
Теперь из  (5),(6) и $I\otimes I\in W_{m,n}$ получим
$${T^p}\otimes {T^p}\in W_{m,n}, \ p\in\Z.$$
Положим $p=n$, тогда из инвариантности $W_{m,n}$  относительно 
 оператора $T^m\otimes T^n$  приходим к $$T^{n-m}\otimes I\in W_{m,n}.$$
Из (4) в силу результатов \cite{Ry} вытекает, что все ненулевые степени  $T^p$
имеют простой спектр.   Рассмотрим циклический вектор $f$  для $T$, он 
является циклическим для всех степеней   $T^p$.
Поэтому из принадлежности   
$$T^{k(n-m)}f\otimes f\in  C_{f\otimes f}$$
для всех $k\in\Z$, где $C_{f\otimes f}$ -- циклическое пространство 
с циклическим вектором ${f\otimes f}$
оператора $T^m\otimes T^n$,   следует включение 
$$ L_2\otimes f \subset  C_{f\otimes f}.$$
Так как $ C_{f\otimes f}$  инвариантно относительно  $T^m\otimes T^n$,    имеем
  $$L_2\otimes T^{kn}f \subset  C_{f\otimes f}$$
для всех $k\in\Z$.
Но  $f$ -- циклический вектор для оператора $T^n$, что дает нужное равенство
 $$L_2\otimes L_2 =  C_{f\otimes f}.$$
Простота спектра $T^m\otimes T^n$ установлена.

Из известных  результатов 
(см. \cite{Ry},\cite{A})  вытекает,  что $T^m\otimes T^m$, $m\neq 0$,  имеет однородный  
сингулярный спектр кратности 2.

Осталось отметить, что диссипативность произведений  $T\times T^n$  при  $n> N$ 
вытекает очевидным образом из условий
$$(A\times B)\cap  (T\times T^n)^k (A\times B)=\phi,$$  
которые выполнены 
для любых  $\xi_j$-измеримых  множеств  $A,B$ (конечной меры)  и  всех   $k>h_j$. 
Теорема доказана.  

Похожие рассуждения приводят к следующему результату.
\vspace{4mm}

\bf  Теорема 3. \it  Для любого $M\subset \N $ найдется  бесконечный автоморфизма $T$ такой, что  
  спектр произведения $T\otimes T^m$     простой  сингулярный  при $m\in M$ и
     счетнократный лебеговский  при $m\notin M$.\rm
\vspace{2mm}

\rm

\section{Заключительные замечания.} 

 Для преобразования $T$ из \cite{KuR},  верно, что  произведения 
$T\otimes T^{n}$  имеют как лебеговский, так и сингулярный спектр для бесконечных множеств
 значений  $n$.  Похожими методами    для некоторого бесконечного потока $T_t$ ранга 1  (см.  \cite{KR})  можно
получить необычное  поведение   
 спектральных кратностей  произведений $T_t\otimes T_{at}$, что будет контрастировать 
 с результатами работ \cite{LR},\cite{K}. 

Действия преобразований, фигурировавших в доказательстве теоремы 2, обладают  
 слабым замыканием с нетривиальной структурой.
   Последовательность $n(k)$, для которой  $T^{n(k)}$ имеет
 ненулевой предел,  обязана иметь вид 
$$n(k)= s+   \alpha_1 h_{j_1(k)} + \alpha_2 h_{j_2(k)}\dots +  \alpha_m h_{j_m(k)},$$ 
где $m,s$ -- некоторые константы, а  $\alpha_i$ принимают целые значения от   $-N$ до $N$. 
Соответствующие слабые пределы являются полиномами. Из теоремы 1 вытекает, что они обладают
минимальными самоприсоединениями, в частности, тривиальным централизатором.

Эффекты, обнаруженные для бесконечных преобразований, в модифицированном виде 
наследуются гауссовскими и пуассоновскими  надстройками, действующими на вероятностном пространстве. 
В качестве примера приведем следующее  утверждение.
\vspace{2mm}

\bf  Теорема 4. \it  Для любого $M\subset \N $ найдется  пуассоновская надстройка $P$ такая, что  
  спектр произведения $P\otimes P^m$       сингулярный  при $m\in M$,  а при $m\notin M$
    в спектре имеется  счетнократная лебеговская компонента.\rm
\vspace{2mm}

Доказательство этой теоремы и ее вариаций будут даны в отдельной работе.
\vspace{2mm}

Тематика этой статьи  обсуждалась на  семинарах, работающих  под руководством 
Б.М. Гуревича, А.А. Давыдова, Ю.А. Неретина, В.И. Оселедца, \\ С.А. Пирогова, Я.Г. Синая, А.М. Степина.
Автор благодарен им за интерес к его работе.

\end{document}